\documentclass[reprint,aps,prx,floatfix,superscriptaddress]{revtex4-2-shared-email}
\usepackage{amsmath}
\usepackage{amssymb}
\usepackage{amsfonts}
\usepackage{physics}
\usepackage{graphicx}
\usepackage{algorithm}
\usepackage{algpseudocode}
\usepackage[breaklinks = true,
            colorlinks = true,
            linkcolor  = blue,
            citecolor  = blue,
            urlcolor   = blue]
            {hyperref}

\usepackage[capitalise,nameinlink]{cleveref}
\crefname{figure}{Fig.}{Figs.}
\Crefname{figure}{Figure}{Figures}
\crefname{table}{Table}{Tables}
\Crefname{table}{Table}{Tables}
\crefname{section}{Sec.}{Secs.}
\Crefname{section}{Section}{Sections}
\crefname{subsection}{Sec.}{Secs.}
\Crefname{subsection}{Section}{Sections}
\crefname{equation}{Eq.}{Eqs.}
\Crefname{equation}{Equation}{Equations}
\crefname{appendix}{Appendix}{Appendices}
\Crefname{appendix}{Appendix}{Appendices}

\usepackage{etoolbox}
\apptocmd{\sloppy}{\hbadness 10000\relax}{}{}

\begin{document}
\title{A Polynomial-Scaling PDE Solver with Entanglement-Basis Tensor Networks} 

\author{Abhijatmedhi Chotrattanapituk}
\thanks{These authors contributed equally}
\email{c\_earth@mit.edu}
\affiliation{Quantum Measurement Group, MIT, Cambridge, Massachusetts 02139, USA}
\affiliation{Department of Electrical Engineering and Computer Science, MIT, Cambridge, Massachusetts 02139, USA}

\author{Michael J. Landry}
\thanks{These authors contributed equally}
\affiliation{Quantum Measurement Group, MIT, Cambridge, Massachusetts 02139, USA}
\affiliation{Department of Nuclear Science and Engineering, MIT, Cambridge, Massachusetts 02139, USA}
\affiliation{Department of Physics, MIT, Cambridge, MA 02139, USA}

\author{Chu-Liang Fu}
\affiliation{Quantum Measurement Group, MIT, Cambridge, Massachusetts 02139, USA}
\affiliation{Department of Nuclear Science and Engineering, MIT, Cambridge, Massachusetts 02139, USA}

\author{Mingda Li}
\email{mingda@mit.edu}
\affiliation{Quantum Measurement Group, MIT, Cambridge, Massachusetts 02139, USA}
\affiliation{Department of Nuclear Science and Engineering, MIT, Cambridge, Massachusetts 02139, USA}

\begin{abstract}
We develop a finite element method (FEM) for partial differential equation (PDE) solver based on the entanglement-basis representation introduced in our companion work. By lifting non-linear finite-element equations into an augmented coefficient space, the governing PDE together with boundary, initial, and inter-element constraints can be expressed through a unified quadratic residual minimization. Although this augmented space grows exponentially with the number of elements, its tensor-product structure allows it to be represented efficiently using tensor networks. Using the matrix product state (MPS) as a concrete example, we show that density matrix renormalization group (DMRG) sweeps enable element-by-element optimization without explicitly constructing the full augmented space. For bounded bond dimension, the resulting computational cost scales polynomially with the number of finite elements. We extend the framework to time-dependent problems through implicit temporal discretization and demonstrate convergence under both mesh and polynomial refinement using diffusion equations.
\end{abstract}

\maketitle

\section{Introduction}\label{sec:introduction}
Partial differential equations (PDEs) underpin quantitative modeling across physics and engineering. Because analytical solutions are available only in restricted cases, practical applications rely on numerical methods such as finite difference (FD)~\cite{leveque2007finite, strikwerda2004finite, thomas1995numerical}, finite element (FE)~\cite{ciarlet2002finite, brenner2008mathematical, hughes2000finite, zienkiewicz2013finite}, and spectral element (SE) methods~\cite{patera1984spectral}. These approaches work by replacing any continuous PDE by a finite algebraic problem via different discretization schemes. The methods also provide systematic convergence under mesh or basis refinement causing their computational cost to grow with resolution, dimensionality, and the complexity of the governing PDE.

Tensor networks (TNs) provide a complementary framework for representing correlated systems efficiently~\cite{orus_2014}. When correlations among the degrees of freedom are sufficiently compressible, such as matrix product states (MPS)~\cite{Fannes_1992,Klumper_1991,Klumper_1993}, projected entangled pair states (PEPS)~\cite{verstraete_2004_A}, and multiscale entanglement renormalization ansatz (MERA)~\cite{Vidal_2007_B}, a TN can approximate a high dimensional state using far fewer parameters than a direct representation. This structure also permits the use of established variational algorithms, including the density matrix renormalization group (DMRG)~\cite{white_1992,white_1993}, time-evolving block decimation (TEBD)~\cite{Vidal_2004}, regularized time evolution, and tree and multi-grid methods\cite{Lubasch_2018,shi2006classical}, for optimizing correlated high-dimensional systems. Therefore, many of recent works have explored TNs as compressed representations for PDEs' numerical solutions~\cite{Bachmayr_2016,Gourianov_2022,Ye_2022,Garca_Ripoll_2021}.

In our companion work, \emph{Tensor-Network Finite Elements for Analytic Operator Equations}, we introduced a unified TN-FE formulation for analytic operator equations (OEs)\cite{chotrattanapituk2026tensor} which is a generalize equation that covers differential, integral, non-linear, memory, and delay operators. The construction maps a discretized operator equation to a hierarchy of multi-linear interaction tensors and represents the resulting coefficient state as a TN. This places broad classes of OEs within a common algebraic framework, incorporates correlations between FE directly into the solution representation, and connects FE discretization with classical and quantum-inspired TN algorithms. In contrast to the companion work, which establishes the general OE framework, the present paper focuses on the explicit implementation and numerical consequences of this formulation. 

The general OE formalism in our companion work does not, however, resolve several practical choices that arise when actually solving PDEs. In conventional FEM, the approximation space, weak formulation, continuity requirements, boundary terms, and discrete operators are commonly designed in relation to the differential order and structure of the specific PDEs\cite{brenner2008mathematical}. Conforming methods build inter-element compatibility into the global basis, whereas discontinuous formulations recover it through fluxes, penalties, or stabilization terms\cite{arnold2002unified}. Consequently, although FEM provide a general discretization philosophy, their concrete realization often remains problem dependent.

In this work, we reorganize these requirements within the TN-FE framework. We begin with independent local representations on each element and impose the governing PDE, inter-element continuity, and boundary conditions as components of variational objective. Continuity and boundary constraints are therefore introduced through explicit residual operators rather than being built entirely into a globally conforming basis. This separates the choice of local approximation functions from the conditions that couple elements and constrain the boundary. The resulting formulation retains the systematic local approximation of FEM while allowing the same TN solver to be used with different elemental representations and boundary relations. Analytic nonlinearities are treated through the augmented coefficient construction introduced in the companion work, enabling their inclusion seamlessly within the same variational framework. 

The remainder of this paper is organized as follows. In \cref{sec:notations}, we introduce the mathematical notations used throughout this works. In \cref{sec:tnfe_overview}, we summarize important results of TN-FE framework from our companion work that will be used without elaborate proof. In \cref{sec:pde_discretize}, we apply TN-FE discretization scheme on a general PDE with some concrete basis choices. We, then, add on inter-element continuities and boundary conditions as additional terms in the residual target optimization in \cref{sec:unify_residual}. In \cref{sec:target_optimization}, we further specialize to one-dimensional for practical demonstration with MPS ansatz and DMRG algorithm. We perform the complexity analysis in \cref{sec:complexity}. In \cref{sec:demonstration}, we demonstrate the method performance through a concrete PDE problem. Finally, in \cref{sec:conclusion}, we discuss implications, limitations, and future perspective of the proposed approach.

\section{Notations}\label{sec:notations}
For the ease of adopting results from our companion work, we will utilize the same mathematical notations as following. First, any non-scalar tensor will be represented with bold face Romans, e.g., $\mathbf{a}$, while scalar and tensor components will be represented with italic Romans, e.g., $a_{\alpha\beta}^{\gamma}$ where Greek characters represent tensor indices. If any index indicates labeling instead of tensor components upright Romans will be used instead, e.g., $a_\mathrm{ij}$.

We use bracket notation to represent an array of the expression inside that got iterated over each of specified indices, i.e., if $\alpha\in\mathbb{Z}_+$, then $\left[a_{\alpha\beta}\right]_{\alpha} = \left[a_{1\beta}, a_{2\beta}, \cdots\right]$ where the subscript $\alpha\gamma$ can be dropped if it is clear by context. We also overload the bracket notation to the indices as well, e.g., $\mathbf{a} = \left[a^{\alpha_1\alpha_2\ldots}\right]=\left[a^{\left[\alpha_\mathrm{i}\right]}\right]$. 

We also utilize Einstein summation for repeated indices. Finally, we use ``$\otimes$'' as tensor product, ``$\odot$'' as tensor contraction, ``$\oplus$'' as tensor space direct sum, ``$\boxplus$'' as tensor block direct sum, and any subscript indices on these symbols indicate that the operations only act on that specific indices.

\section{Overview of the Tensor-Network Finite Element  Framework}
\label{sec:tnfe_overview}

\subsection{Finite Element Representation of Analytic Operator Equations}
Let $\Omega$ denotes the domain of independent variable $\mathbf{q}$, with $\mathcal{H}(\Omega)$ denotes the admissible solution space, and $\mathcal{R}(\Omega)$ denote the residual space. A general OE over $\Omega$ is
\begin{equation}
    \mathbf{L}[\mathbf{u}] = \mathbf{0}
\end{equation}
where $\mathbf{L}:\mathcal{H}(\Omega)\rightarrow\mathcal{R}(\Omega)$ and $\mathbf{u}=\mathbf{u}(\mathbf{q})\in\mathcal{H}(\Omega)$. If we partition $\Omega$ into $N$ elements $\Omega_\mathrm{n}$ for $\mathrm{n}\in\{1,\ldots,N\}$ and introduce local basis functions $b_{\mathrm{n},\beta_\mathrm{n}}$ supported on $\Omega_\mathrm{n}$ with assumption that all elements contain the same number of basis functions, the element-wise approximation is
\begin{equation}
    \mathbf{u} \approx \mathbf{c}^{\nu\beta} b_{\nu\beta}
\end{equation}
where $\mathrm{n}$ is substituted with a tensor index $\nu$ and $\left[\mathbf{c}^{\nu\beta}\right]$ is expansion coefficient tensor. In a Galerkin weak formulation, the residual is projected onto this test space,
\begin{equation}
    \int_\Omega b_{\nu'\beta'}^\ast\,
    \mathbf{L}\left[\mathbf{c}^{\nu\beta} b_{\nu\beta}\right]\,\mathrm{d}\mathbf{q} = \mathbf{0}
\end{equation}
for all test functions $b_{\nu'\beta'}$. This reduces the OE to a finite algebraic problem for the coefficient tensor.

For an analytic operator, $\mathbf{L}$ may be expanded as a converging Fr\'echet-Taylor series~\cite{boyd1984analytical,boyd1985fading,schetzen1980volterra,mujica1986complex,dineen1999complex} of the form
\begin{equation}\label{eq:aoe}
    \mathbf{L}[\mathbf{u}]
    = \sum_{\mathrm{i}=0}^\infty\dfrac{1}{\mathrm{i}!}\int_{\Omega^\mathrm{i}}\mathbf{k}_{\mathrm{i}}^{\left[\eta_\mathrm{j}\right]}\odot\bigotimes_{\mathrm{j}=1}^\mathrm{i}\partial_{\mathbf{q}_\mathrm{j}}^{\eta_\mathrm{j}}\mathbf{u}_\mathrm{j}\,\mathrm{d}\mathbf{q}_\mathrm{j}
\end{equation}
where $\eta_\mathrm{j}\in\mathbb{N}_0$, $\mathbf{u}_\mathrm{j} = \mathbf{u}(\mathbf{q}_\mathrm{j})$, and $\mathbf{k}_{\mathrm{i}}^{\left[\eta_\mathrm{j}\right]}$ is an $\mathrm{i}$-linear interaction kernel with proper contraction to each $\partial_{\mathbf{q}_\mathrm{j}}^{\eta_\mathrm{j}}\mathbf{u}_\mathrm{j}$. After inserting this, the weak form becomes a hierarchy of multi-linear algebraic contributions,
\begin{equation}\label{eq:wf_fe_aoe}
    \sum_{\mathrm{i}} \mathbf{g}_{\mathrm{i},\nu'\beta', \left[\nu_\mathrm{j}\beta_\mathrm{j}\right]}\odot\bigotimes_{\mathrm{j}=1}^\mathrm{i}\mathbf{c}^{\nu_\mathrm{j}\beta_\mathrm{j}}
    = \mathbf{0}
\end{equation}
where $\left[\mathbf{g}_{\mathrm{i},\nu'\beta', \left[\nu_\mathrm{j}\beta_\mathrm{j}\right]}\right]_{\nu'\beta', \left[\nu_\mathrm{j}\beta_\mathrm{j}\right]}$ is the basis projection tensor of the $\mathrm{i}^\text{th}$ order correlation term of $\mathbf{L}$.

\subsection{Augmented coefficient space}
In conventional finite elements, the global approximation space is assembled as a direct sum of the local spaces:
\begin{equation}
    \Phi = \bigoplus_{\nu=1}^N\Phi_\nu
\end{equation}
where $\Phi_\nu = \text{Span}\left(\left\{b_{\nu\beta}\mid \beta\in\{1,\ldots,B\}\right\}\right)$ and $B$ is the number of basis functions in each element. In this representation space, any linear $\mathbf{L}$ can be represented directly as a matrix equation. However, for non-linear analytic $\mathbf{L}$, additional intelligent variable substitution and/or iterative root finding algorithm is required.

The TN-FE construction can represent these correlation linearly at the cost of enlarging the representation space by lifting $\Phi$ to its Fock space defined as
\begin{equation}\label{eq:fock}
    \Psi
    = \bigoplus_{\mathrm{i}=0}^\infty \Phi^{\otimes\mathrm{i}}\,.
\end{equation}
With this construction, any product $\bigotimes_{\mathrm{j}}\mathbf{c}^{\nu_\mathrm{j}\beta_\mathrm{j}}$ can be replaced by an independent coefficient in $\Psi$. Furthermore, with the multi-linear nature of $\mathbf{L}$, we have gauge freedom to permute terms in the product such that they are in the order of element label. This means that by introducing, for each element, empty basis function $b_{\nu0}$ and letting $J$ be the largest correlation order of $\mathbf{L}$, we can work in the basis system that, for each element $\nu$, $B+1$ basis functions (including $b_{\nu0}$) can be allocated to $J$ slots. Hence, \cref{eq:wf_fe_aoe} can be replaced with
\begin{equation}\label{eq:wf_aug_fe_aoe}
    \mathbf{G}_{\left[\beta'_{\nu\mathrm{j}}\right],\left[\beta_{\nu\mathrm{j}}\right]}\odot\mathbf{c}^{\left[\beta_{\nu\mathrm{j}}\right]} = \mathbf{0}
\end{equation}
where $\mathrm{j}$ indicates the slot number and $\left[\mathbf{G}_{\left[\beta'_{\nu\mathrm{j}}\right],\left[\beta_{\nu\mathrm{j}}\right]}\right]$ is the projection tensor of $\mathbf{L}$ in this augmented space. For cleanliness, we will directly use $\mathbf{c}$ and $\mathbf{G}$ for coefficient and operator projection tensor, respectively, without the indices and brackets.

\begin{figure}[!htb]
    \centering
    \includegraphics[width=\columnwidth]{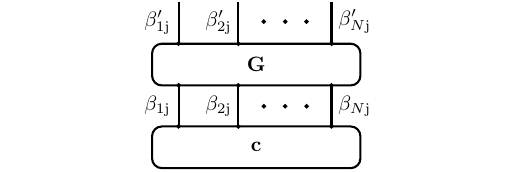}
    \caption{TN representing the residual term of the linear problem constructed from a general analytic OE with TN-FE framework (\cref{eq:wf_aug_fe_aoe}). Each box node represent a tensor in the expression. Free indices are represented with hanging bonds while index contractions are represented with connecting bonds. Each bond can correspond to a group of indices. In this case, each bond represents $J$ indices of the same element. (Adapted from \cite{chotrattanapituk2026tensor})}
    \label{fig:residual}
\end{figure}

These results already capture the idea behind TN-FE formulation, but introducing TN visualization can make the application of the method more practical and intuitive. For instance, the left hand side of \cref{eq:wf_aug_fe_aoe} can be neatly represented in TN format as shown in \cref{fig:residual} where the box nodes representing $\mathbf{G}$ and $\mathbf{c}$ are connected through the bonds that represent index contraction between them.

\subsection{Tensor Network Compression}
From the construction of Fock space (\cref{eq:fock}), direct storage of $\mathbf{c}$ scales exponentially with number of elements and correlation order as $\mathcal{O}(B^{NJ})$. Hence, the Fock space can only be used as a formal representation for correlated system to be further compressed by any TN ansatz.

\begin{figure}[!htb]
    \centering
    \includegraphics[width=\columnwidth]{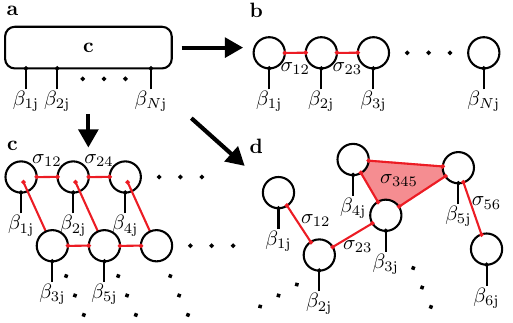}
    \caption{TN with locality ansatzes of coefficient tensor under different geometric restriction. Here, $\sigma$'s indicate tensor indices responsible for each local correlation. \textbf{(a)} No locality assumption: same coefficient tensor as in \cref{fig:residual}. \textbf{(b)} Matrix product state (MPS): elements align in one-dimensional chain with correlation only between neighbors. \textbf{(c)} Projected entangled pair state (PEPS):elements align in two-dimensional grid with correlation between adjacent nodes along the grid lines. \textbf{(d)} Higher-order TN assumption where shaded regions represent index hyper-contraction. (Adapted from \cite{chotrattanapituk2026tensor})}
    \label{fig:tn_ansatz}
\end{figure}

There is no definite rule to design a TN ansatz, but one popular assumption is the locality of correlation where the correlation is expected to only appear between nearby elements. \Cref{fig:tn_ansatz} illustrates representative TN ansatzes of $\mathbf{c}$ (\cref{fig:tn_ansatz}a) with locality assumption for different geometric restriction: matrix product state (MPS, \cref{fig:tn_ansatz}b) for one-dimensional element chain with only correlation between neighbors, projected entangled pair state (PEPS, \cref{fig:tn_ansatz}c) for two-dimensional grid with only correlation between axial neighbors, and higher-order TN (\cref{fig:tn_ansatz}d) for more than pair-wise correlation.

If the maximum bond dimension (range of $\sigma$'s) in a TN is $S$, and the maximum bond degree (including hyper-contraction) in the TN is $M$, then the storage requirement is 
\begin{equation}
    \mathcal{O}(NB^JS^M)
\end{equation}
instead of $\mathcal{O}(B^{NJ})$, i.e., not exponentially growing with $N$. More generally, the storage remains polynomial in $N$ when the local correlation order $J$, bond dimension $S$, and network degree $M$ remain bounded.

\subsection{Variational solution of the discretized equation}
For a coefficient tensor $\mathbf{c}$, define the projected residual
\begin{equation}\label{eq:aoe_res}
    \mathbf{h}_{\left[\beta'_{\nu\mathrm{j}}\right]}(\mathbf{c}) = \mathbf{G}_{\left[\beta'_{\nu\mathrm{j}}\right],\left[\beta_{\nu\mathrm{j}}\right]}\odot\mathbf{c}^{\left[\beta_{\nu\mathrm{j}}\right]}\,.
\end{equation}
The TN-FE solution is obtained by minimizing the weighted residual norm
\begin{equation}\label{eq:aoe_wrn}
    J(\mathbf{c})
    = W^{\left[\beta'_{\nu\mathrm{j}}\right],\left[\beta_{\nu\mathrm{j}}\right]}\left(\mathbf{h}^\dagger_{\left[\beta'_{\nu\mathrm{j}}\right]}(\mathbf{c})\odot \mathbf{h}_{\left[\beta_{\nu\mathrm{j}}\right]}(\mathbf{c})\right)\,,
\end{equation}
where $W^{\left[\beta'_{\nu\mathrm{j}}\right],\left[\beta_{\nu\mathrm{j}}\right]}$ is a positive-definite metric on the residual space. Since we parametrize the search space to that of particular TN ansatz manifold $\mathcal{T}_{\mathrm{TN}}$, i.e., 
\begin{equation}
    \mathbf{c}_{\mathrm{TN}}
    =
    \arg\min_{\mathbf{c}\in\mathcal{T}_{\mathrm{TN}}}
    J(\mathbf{c})\,,
\end{equation}
the optimization generally become non-convex. However, if we fix all but one node of $\mathbf{c}$'s TN ansatz, the single-site update reduces to a convex weighted least-squares problem again. This is part of the reason for the success of the alternate optimization algorithm like density matrix renormalization group (DMRG) for MPS ansatz.

\subsection{Specialization to PDE}
The formulation above is for a general analytic OEs, but PDE is actually a sub-class of OE. Hence, it is possible to reduce any analytic PDE problems into OE problems and utilize TN-FE formulation. In particular, a PDE can be written as
\begin{equation}
    \mathbf{g}\left(\mathbf{q},\left[\partial^\eta_{\mathbf{q}}\mathbf{u}\right]\right)
    =
    \mathbf{0}\,.
\end{equation}
If $\mathbf{g}$ is analytic, then it can be expanded as a convergent series 
\begin{equation}\label{eq:reduce_PDE}
    \mathbf{g}\left(\mathbf{q},\left[\partial^\eta_{\mathbf{q}}\mathbf{u}\right]\right)
    = \sum_{\mathrm{i}}\dfrac{1}{\mathrm{i}!}\mathbf{a}_{\mathrm{i}}^{\left[\eta_\mathrm{j}\right]}\odot\bigotimes_{\mathrm{j}=1}^\mathrm{i}\partial_{\mathbf{q}}^{\eta_\mathrm{j}}\mathbf{u}\,,
\end{equation}
which is a special case of \cref{eq:aoe} when
\begin{equation}
    \mathbf{k}_{\mathrm{i}}^{\left[\eta_\mathrm{j}\right]} = \mathbf{a}_{\mathrm{i}}^{\left[\eta_\mathrm{j}\right]}\bigotimes_{\mathrm{j}=1}^{\mathrm{i}}\delta(\mathbf{q}_\mathrm{j}-\mathbf{q})\,.
\end{equation}
The new emphasis of the present work is the explicit construction and implementation of this PDE solver.

\section{PDE Discretization}\label{sec:pde_discretize}

\subsection{Target Equation}\label{sec:target_eq}
In this work, we are considering PDEs that has temporal dependence and, for convenience, we explicitly single out the time variable $t$ as
\begin{equation}
    \mathbf{g}\left(t, \mathbf{q}, \left[\partial^\eta_{\mathbf{q}}\partial^\zeta_{t}\mathbf{u}\right]\right) = \mathbf{0}\,.
\end{equation}
There is nothing wrong with treating $t$ as a dimension of $\mathbf{q}$. We decided to treat temporal discretization differently from the spatial part to be able to utilize MPS ansatz and DMRG algorithm on the one-dimensional initial value problem (IVP) demonstration later. If $\mathbf{g}$ is also analytic, we can rewrite the PDE as 
\begin{multline}\label{eq:td_apde}
    \mathbf{g}\left(t, \mathbf{q}, \left[\partial^\eta_{\mathbf{q}}\partial^\zeta_{t}\mathbf{u}\right]\right)\\
    = \sum_{\mathrm{i}}\dfrac{1}{\mathrm{i}!}\mathbf{a}_{\mathrm{i}}^{\left[\zeta_\mathrm{k}\right], \left[\eta_\mathrm{j}\right]} \odot\bigotimes_{\mathrm{j}=1}^\mathrm{i}\bigotimes_{\mathrm{k}=\mathrm{j}+1}^\mathrm{i}\partial_{\mathbf{q}}^{\eta_\mathrm{j}}\partial_{t}^{\zeta_\mathrm{k}}\mathbf{u}\,.
\end{multline}

\subsection{Temporal discretization}
For an IVP, given $\widetilde{\mathbf{u}}$ and some of its time derivatives at time $t$, we want to determine $\mathbf{u}$ and those time derivatives at time $t+\Delta$ where $\Delta$ is time step size. The most common way is to apply finite difference method (FDM) to the temporal part~\cite{leveque2007finite}. One can choose to use any FDM scheme of their preference with varying degree of convergence criteria and computational costs, but the core idea of the following derivation remains the same. 

For simplicity, we decided to use the implicit Euler's scheme. In other words, \cref{eq:td_apde} is replaced with
\begin{multline}
    \mathbf{0} = \sum_{\mathrm{i}}\dfrac{1}{\mathrm{i}!}\mathbf{a}_{\mathrm{i}}^{\left[\zeta_\mathrm{k}\right],\left[\eta_\mathrm{j}\right]}\odot\bigotimes_{\mathrm{j}=1}^\mathrm{i}\bigotimes_{\mathrm{k}=\mathrm{j}+1}^\mathrm{i}\\
    \partial_{\mathbf{q}}^{\eta_\mathrm{j}}
    \begin{cases}
        \dfrac{1}{\Delta}\left(\partial_{t}^{\zeta_\mathrm{k}-1}\mathbf{u}-\partial_{t}^{\zeta_\mathrm{k}-1}\widetilde{\mathbf{u}}\right)&, \zeta_\mathrm{k}\neq 0\\
        \mathbf{u}&, \zeta_\mathrm{k} = 0
    \end{cases}\,.
\end{multline}
With some term arrangement, the equation can be written as a time-dependent PDE on $\left[\partial_{t}^{\zeta}\mathbf{u}\right]$, i.e., it is always possible to simplify the equation to only contain up to first temporal derivatives.

With FDM scheme, we will apply the method iteratively to solve for solution at various time steps
\begin{equation}
    t^{(\mathrm{t})} = t + \mathrm{t}\Delta
\end{equation}
where $\mathrm{t}\in\mathbb{Z}_{0+}$. With this new notation introduction, our target equation is
\begin{multline}\label{eq:ivp_form}
    \mathbf{0} = \sum_{\mathrm{i}}\dfrac{1}{\mathrm{i}!}\mathbf{a}_{\mathrm{i}}^{\left[\zeta_\mathrm{k}\right],\left[\eta_\mathrm{j}\right]}\odot\bigotimes_{\mathrm{j}=1}^\mathrm{i}\bigotimes_{\mathrm{k}=\mathrm{j}+1}^\mathrm{i}\\
    \partial_{\mathbf{q}}^{\eta_\mathrm{j}}
    \begin{cases}
        \dfrac{\mathbf{u}^{(\mathrm{t}+1)}-\mathbf{u}^{(\mathrm{t})}}{\Delta}&, \zeta_\mathrm{k} = 1\\
        \mathbf{u}^{(\mathrm{t}+1)}&, \zeta_\mathrm{k} = 0
    \end{cases}\,.
\end{multline}

We can also see that $\mathbf{u}^{(\mathrm{t})}$ and $\Delta$ can be absorbed to the term with lower correlation order and interaction kernel, respectively. This basically recover the analytic PDE, i.e., our method on IVP can be applied to any PDE.

\subsection{Building Block Operators}
One can directly apply TN-FE formulation on the IVP of interest, but it means that the calculation of $\mathbf{G}$ need to be re-calculated from scratch at every time step because the multi-linear interaction kernel can be time dependent and the coefficient tensor of the previous time step got updated on the next iteration. However, because everything else in the PDE is the same, it should be possible to pre-compute those parts of the PDE into a concise weak form ready for connecting to the time-dependent parts.

\subsubsection{Derivative}
Since spatial derivative is a local operators, if we have well-behaved basis function (can approximate every viable function in the element well), the results of derivative on $\Phi_\mathrm{n}$ can also be well approximated with a member of $\Phi_\mathrm{n}$. In other words, we can approximate spatial derivative with a projection tensor.

Consider an arbitrary spatial derivative of a basis function $b_{\mathrm{n},\beta}$ which can be approximated as
\begin{equation}
    \partial_{\mathbf{q}}b_{\mathrm{n},\beta} \approx m^{\beta'}b_{\mathrm{n},\beta'}\,.
\end{equation}
Applying weak formulation gives an approximated projection tensor of $\partial_{\mathbf{q}}$ on the element $\mathrm{n}$ as
\begin{equation}\label{eq:derivative_projection}
     \mathbf{D}_{\mathrm{n},\beta'\beta}\approx Q^{-1}_{\mathrm{n},\beta'\beta''}\mathbf{P}^{(\text{der})}_{\mathrm{n},\beta''\beta}\,.
\end{equation}
where
\begin{align}
    Q_{\mathrm{n},\beta''\beta'} &= \int_{\Omega_\mathrm{n}} b^\ast_{\mathrm{n},\beta''}b_{\mathrm{n},\beta'}\,\mathrm{d}\mathbf{q}\\
    \mathbf{P}^{(\text{der})}_{\mathrm{n},\beta''\beta} &= \int_{\Omega_\mathrm{n}} b^\ast_{\mathrm{n},\beta''}\partial_{\mathbf{q}}b_{\mathrm{n},\beta}\,\mathrm{d}\mathbf{q}\,.
\end{align}

We can also visualize it in TN form as shown in \cref{fig:diff_mul_fac_bp}a. Note that the label $\mathrm{n}$ is needed since each element has different sets of basis functions. Also, $\mathbf{D}_{\mathrm{n},\beta'\beta}$ is a tensor with components (tensor index $\chi$) correspond to different axis of differentiation.

\subsubsection{Multiplication}
With similar method to derivative, we can achieve the approximation of spatial multiplication with projection tensor 
\begin{equation}\label{eq:multiplication_projection}
    R_{\mathrm{n},\beta'\beta_1\beta_2} \approx Q^{-1}_{\mathrm{n},\beta'\beta''}P^{(\text{mul})}_{\mathrm{n},\beta''\beta_1\beta_2}\,,
\end{equation}
where
\begin{equation}
    P^{(\text{mul})}_{\mathrm{n},\beta''\beta_1\beta_2}=\int_{\Omega_\mathrm{n}} b^\ast_{\mathrm{n}\,\beta''}b_{\mathrm{n},\beta_1}b_{\mathrm{n},\beta_2}\,\mathrm{d}q\,,
\end{equation}
with visual depiction in \cref{fig:diff_mul_fac_bp}b. 

\subsubsection{Constant Factor}
Constant factor can imply scaling of the whole function by a constant factor or tensor contraction to the indices that are not basis indices $\beta$'s, i.e., the action on basis is simply identity projection. If said constant tensor is $\mathbf{f}$, with similar method to derivative and multiplication, the projection tensor for constant factor is
\begin{equation}
    \mathbf{F}_{\mathrm{n},\beta'\beta} = Q^{-1}_{\mathrm{n},\beta'\beta''}Q_{\mathrm{n},\beta''\beta}\mathbf{f} = \delta_{\beta'\beta}\mathbf{f}
\end{equation}
where $\delta$ is Kronecker delta, i.e., $\mathbf{f}$ only acts on the non-basis indices, with visualization in \cref{fig:diff_mul_fac_bp}c. This component might look trivially simple but it is the foundation for applying our framework toward PDE on tensor value functions such as is Electrodynamics.

\subsubsection{Boundary Projection}
\begin{figure}[htb!]
    \centering
    \includegraphics[width=\columnwidth]{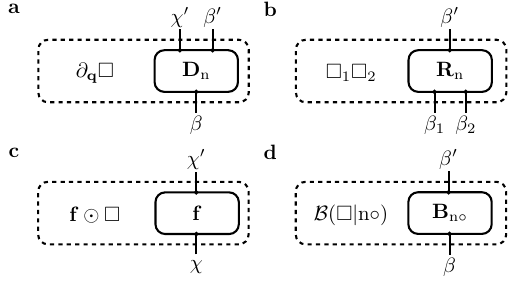}
    \caption{TN representation of building block operations: \textbf{(a)} spatial derivative, \textbf{(b)} multiplication, \textbf{(c)} constant factor, and \textbf{(d)} projection of element $\mathrm{n}$ to boundary with $\circ$. Here, $\beta$'s are basis and $\chi$'s are non-basis indices. ($\mathbf{a}$ and $\mathbf{b}$ are adapted from \cite{chotrattanapituk2026tensor})}
    \label{fig:diff_mul_fac_bp}
\end{figure}

Later, we will integrate boundary conditions, and inter-element continuity into our TN variational target. Hence, it is crucial to look at the projection of elemental solution to element edges. Without lost of generality, let elements $\mathrm{n}$ and $\circ$ (can be environment) be each other neighbor. The boundary between them is $\Omega_{\mathrm{n}\circ}=\Omega_{\mathrm{n}}\cap\Omega_{\circ}$ where we abuse notation for the partitions to include boundaries. Let $\{b_{\mathrm{n}\circ,\beta}\}$ be boundary's basis. Then, the projection tensors of element $\mathrm{n}$ to the boundary with element $\circ$ is
\begin{equation}\label{eq:boundary_projection}
    B_{\mathrm{n}\circ,\beta'\beta} \approx Q^{-1}_{\mathrm{n}\circ,\beta'\beta''}P^{(\text{bp})}_{\mathrm{n}\circ,\beta''\beta}
\end{equation}
where
\begin{align}
    Q_{\mathrm{n}\circ,\beta''\beta'} &= \int_{\Omega_{\mathrm{n}\circ}} b^\ast_{\mathrm{n}\circ,\beta''}b_{\mathrm{n}\circ,\beta'}\,\mathrm{d}\mathbf{q}\\
    P^{(\text{bp})}_{\mathrm{n}\circ,\beta''\beta} &= \int_{\Omega_{\mathrm{n}\circ}} b^\ast_{\mathrm{n}\circ,\beta''}b_{\mathrm{n},\beta}\,\mathrm{d}\mathbf{q}\,.
\end{align}
We can also visualize it in TN form as shown in \cref{fig:diff_mul_fac_bp}d.

\subsubsection{Summation}
There are three places that summation appear in our formalism: the expansion in \cref{eq:aoe}, inter-element continuities (ICs) and boundary conditions (BCs), and the inclusion of said ICs and BCs into the variational target. Among these, the summations can be classified into two categories depending on whether all summands share the same TN structure.

For the cases that all TNs representing the summands have the same structure, we can equivalently represent the summation by creating a new TN with the same structures as summand TNs with each tensor node replaced by the tensor block direct sum of the same nodes from all summand TNs. 

\begin{figure}[htb!]
    \centering
    \includegraphics[width=\columnwidth]{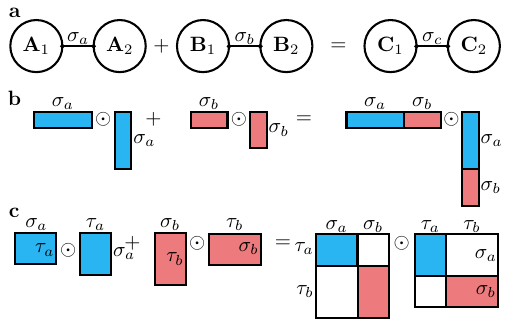}
    \caption{Summation of TNs with the same structure. \textbf{(a)} Diagram representing summation of two TNs with the same structure which can be replace with a single TN with the same structure with each node be the tensor block direct sum of each original node at the same position in the structure, i.e., $\mathbf{C}_\mathrm{i} = \mathbf{A}_\mathrm{i}\boxplus_\sigma\mathbf{B}_\mathrm{i}$. \textbf{(b)} Simple example of \textbf{a} where each node represent a vector. \textbf{(c)} Slightly more complicated example of \textbf{a} where each node represent a matrix and there are two bonding indices completely contract matrices to scalar, i.e., $\mathbf{C}_\mathrm{i} = \mathbf{A}_\mathrm{i}\boxplus_{\sigma\tau}\mathbf{B}_\mathrm{i}$.}
    \label{fig:summation_1}
\end{figure}

To see that this is true, consider a summation between two TNs along a common bond position as shown in \cref{fig:summation_1}a. Since both TNs on the left can be added, they must have the same dimension in all indices except the one that got contracted. We will focus on a particular bond position in the TN shown as $\sigma_a$ and $\sigma_b$, i.e., we contract all remaining bonds in the TN. The summation of the TNs is
\begin{equation}\label{eq:tn_structure_sum}
    \mathbf{A}_{1,\sigma_a}\mathbf{A}_{2,\sigma_a} + \mathbf{B}_{1,\sigma_b}\mathbf{B}_{2,\sigma_b}\,.
\end{equation}

Because both TNs have the same structure and non-bonding dimensions, it is possible to replace the summation indices with $\sigma_c$ which run through all values of $\sigma_a$ followed by values of $\sigma_b$ and construct a new tensor node
\begin{equation}
    \mathbf{C}_{\mathrm{i},\sigma_c} =
    \begin{cases}
        \mathbf{A}_{\mathrm{i},\sigma_c}&, \sigma_c \leq \sigma_a\\
        \mathbf{B}_{\mathrm{i},\sigma_c-\sigma_a}&, \sigma_c > \sigma_a
    \end{cases}
\end{equation}
which is a block direct sum between the summand nodes along bonding index, i.e., $\mathbf{C}_\mathrm{i} = \mathbf{A}_\mathrm{i}\boxplus_\sigma\mathbf{B}_\mathrm{i}$. With this new tensor we can directly see that
\begin{equation}
    \mathbf{C}_{1,\sigma_c}\mathbf{C}_{2,\sigma_c} = \mathbf{A}_{1,\sigma_a}\mathbf{A}_{2,\sigma_a} + \mathbf{B}_{1,\sigma_b}\mathbf{B}_{2,\sigma_b}.
\end{equation}
\Cref{fig:summation_1}b also demonstrates this for the case when all nodes are vectors. It might seems like using space direct sum works here as well, but block version is needed when the bonding has more than one index such as in the case of \cref{fig:summation_1}c where all nodes are matrices that got completely contracted to a scalar. Note that, this analysis also applies when the bond under consideration is connecting two indices of the same tensor.

We can directly apply this results to our formalism. However, in our cases, because all of the summands share the same Fock space coefficient tensor, we can further simplify the summation building block by adding (not block direct sum) some of the tensors of the same node position in the common TN structure. As a concrete example, if the right vectors of both terms on the right of \cref{fig:summation_1}b are the same, we can perform the summation of the left vectors of both terms before contracting to the common right vector. This effectively saves up one tensor index contraction which can cascade into considerable amount of time saved for the more complex cases.

On the other hand, when the underlying TN structure of the summands are not the same, it is more complicated to represent those summations in TN format even if all summands' TN share the same Fock space coefficient tensor. To solve this issue, we will consider the summation between two arbitrary tensors $\mathbf{A}$ and $\mathbf{B}$. To merge them into a single TN, we can consider their tensor product $\mathbf{A}\otimes\mathbf{B}$. Unfortunately, their is no linear tensor operation that can get individual $\mathbf{A}$ or $\mathbf{B}$ out of $\mathbf{A}\otimes\mathbf{B}$. Hence, we will consider an additional augmentation to each of them with block direct sum with 1. By considering their tensor product, we get
\begin{multline}\label{eq:sum_aug_tp}
    (\mathbf{A}\boxplus 1_a)\otimes(\mathbf{B}\boxplus 1_b) = (\mathbf{A}\otimes\mathbf{B}) \oplus (\mathbf{A}\otimes 1_b) \\
    \oplus (1_a\otimes\mathbf{B}) \oplus (1_a\otimes 1_b)
\end{multline}
where the subscripts to $1$'s only indicate the vectors that they are associated to.

\begin{figure}[htb!]
    \centering
    \includegraphics[width=\columnwidth]{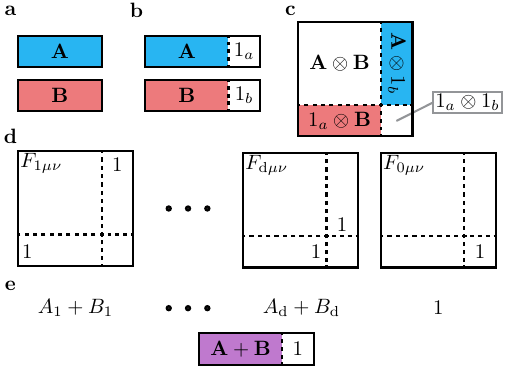}
    \caption{Summation of two arbitrary TNs (not necessarily the same structure) that contracted to vectors $\mathbf{A}$ and $\mathbf{B}$ with $\mathrm{d}$ dimensions. \textbf{(a)} Diagrams representing $\mathbf{A}$ and $\mathbf{B}$. \textbf{(b)} Augmentation of the vectors. \textbf{(c)} Tensor product of the augmented vectors which contains subspaces that are linear to $\mathbf{A}$ and $\mathbf{B}$. \textbf{(d)} Filter $\mathbf{F}$ that extracts and adds components of $\mathbf{A}$ and $\mathbf{B}$. The diagrams show each of the $\mathrm{d}+1$ layers of $\mathbf{F}$ with only non-zero numbers shown explicitly. \textbf{(e)} Extraction by each layer of $\mathbf{F}$ that corresponds to a component of $(\mathbf{A}+\mathbf{B})\boxplus 1$.}
    \label{fig:summation_2}
\end{figure}

We can see that the second and third terms of \cref{eq:sum_aug_tp} are linear to $\mathbf{A}$ and $\mathbf{B}$, respectively. Therefore, we can design a filter $\mathbf{F}$ that extracts and adds them together as
\begin{equation}
    \mathbf{F}\odot\left((\mathbf{A}\boxplus 1_a)\otimes(\mathbf{B}\boxplus 1_b)\right)=(\mathbf{A}+\mathbf{B})\boxplus 1.
\end{equation}
We also illustrated this result in \cref{fig:summation_2} for the case when $\mathbf{A}$ and $\mathbf{B}$ are vectors.

Building from this, consider the case when $\mathbf{A}$ (and $\mathbf{B}$) can be factored out like in \cref{eq:tn_structure_sum}. We can see that
\begin{equation}
    (\mathbf{A}_{1,\sigma_a}\boxplus1_{a1})(\mathbf{A}_{2,\sigma_a}\boxplus1_{a2}) = (\mathbf{A}_{1,\sigma_a}\mathbf{A}_{2,\sigma_a})\boxplus1_a\,,
\end{equation}
i.e., we can apply the augmentation to each node in TN instead of applying after the complete contraction. Here, we also abuse some notation for index $\sigma_a$ which need to run over the index value for $1_{a1}$ and $1_{a2}$ as well. Furthermore, similar to the same structure case, if parts of the summands are the same such as sharing the same Fock space coefficient tensor, we can factor out that parts. For example, if $\sigma_a$ and $\sigma_b$ share the same index range and $\mathbf{A}_{2,\sigma}=\mathbf{B}_{2,\sigma}=\mathbf{E}_\sigma$,
\begin{align}\label{eq:sum_factor}
    &(\mathbf{A}_{1,\sigma_a}\mathbf{A}_{2,\sigma_a} + \mathbf{B}_{1,\sigma_b}\mathbf{B}_{2,\sigma_b})\boxplus1 \nonumber\\
    &= ((\mathbf{A}_{1,\sigma} + \mathbf{B}_{1,\sigma})\mathbf{E}_\sigma)\boxplus 1 \nonumber\\
    &= ((\mathbf{A}_{1,\sigma} + \mathbf{B}_{1,\sigma})\boxplus 1)(\mathbf{E}_\sigma\boxplus 1) \nonumber\\
    &= (\mathbf{F}\odot\left((\mathbf{A}_{1,\sigma}\boxplus 1_{a1})\otimes(\mathbf{B}_{1,\sigma}\boxplus 1_{b1})\right))(\mathbf{E}_\sigma\boxplus 1).
\end{align}
This means that we can isolate the summation locally to only the nodes that are different.

\begin{figure}[htb!]
    \centering
    \includegraphics[width=\columnwidth]{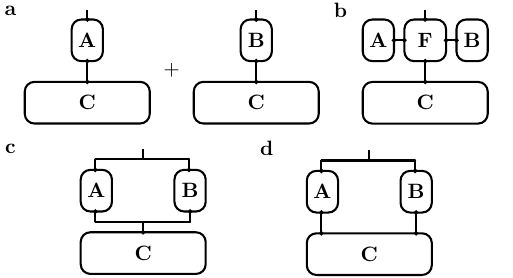}
    \caption{TN representation of summation building block. \textbf{(a)} Summation of two TNs that has a common part $\mathbf{C}$ and the different parts $\mathbf{A}$ and $\mathbf{B}$. \textbf{(b)} Building block representation of $\mathbf{a}$. As shown in \cref{eq:sum_factor}, the common part can be factor out reducing the task to be the summation of the different parts. Here, $\mathbf{F}$ represents filter contracting to $\mathbf{A}$ and $\mathbf{B}$. The block direct sum augmentation is implicitly assumed for cleanliness. \textbf{(c)} Specialization for the cases where $\mathbf{A}$ and $\mathbf{B}$ contract to the same indices of $\mathbf{C}$ where the filter is assumed in the branching and merging of contracting edges. \textbf{(d)} Specialization for the cases where $\mathbf{A}$ and $\mathbf{B}$ contract to the different indices of $\mathbf{C}$.}
    \label{fig:sum}
\end{figure}

Also, the case of summation between TNs with the same structure can be considered a special case of this with a short cut of directly adding each tensor node directly. This technique of applying filter on the tensor product of augmented TNs can be expanded for any other linear operation between TNs. For cleanliness, we will compress out the augmentation from our TN and use simple visualization like in \cref{fig:sum}a-b to represent the summation of two TNs with $\mathbf{C}$ as their common parts while $\mathbf{A}$ and $\mathbf{B}$ are their non-common parts.

Depending on how the summation appear, the filter $\mathbf{F}$ can take multiple form. For our formalism, we specialize the summation for the case that $\mathbf{A}$ and $\mathbf{B}$ are acting on the same or different indices of $\mathbf{C}$ as represented in \cref{fig:sum}c-d, respectively. For further cleanliness, we implicitly suppressed the summation filter for each branching or merging of the contraction edges.

\subsubsection{Compound Network}
With these components, we can construct more complicated IVP operators, e.g., higher-order derivatives, other tensor product variations, and tensor calculus operators. Of course, directly applying the weak formalism for the whole operator $\mathbf{L}$ would give a better approximation at the cost of additional overhead calculations per time step. 

\subsection{Local Basis}
\subsubsection{Independent Local Approximation Spaces}
In standard Galerkin weak formulation for FEM, the elemental basis functions must be chosen such that their assembly form a globally conforming basis. This can be either choosing the basis that is already a globally conforming basis, or impose strict condition on coefficient tensor. For example, if the solution is expected to be continuous, the basis functions for each element must be 0 along all of the element boundaries or the coefficient tensor for that element must be restricted by neighboring elements' coefficient tensors to make sure that the approximated solution is continuous.

Since, in TN-FE formulation, the residual is optimized via quadratic target scheme, boundary conditions, continuity and any required derivative compatibility can be included through additional residual terms. Consequently, the basis functions need not be assembled into a globally conforming basis, and coefficient tensors can be treated as independent degrees of freedom. Here, we will explore two bases that would be of interest since they can have their close form for derivative, multiplication and boundary projection tensors. The constant factor does not involve basis to begin with.

\subsubsection{Fourier Basis}
Consider a one-dimensional domain element $\Omega_\mathrm{n}$ with local coordinate origin $q_\mathrm{n}$, i.e., $q'=q-q_\mathrm{n}$, and fundamental wavelength $l_\mathrm{n}$. We define Fourier basis functions with $\Omega_\mathrm{n}$ support as
\begin{equation}
    b^{(\text{FR})}_{\mathrm{n},\beta} =
    \begin{cases}
    \sin\left(\beta k_\mathrm{n}q'\right)&; \beta > 0\\
    \cos\left(\beta k_\mathrm{n}q'\right)&; \beta \leq 0
    \end{cases}
\end{equation}
where $k_\mathrm{n} = 2\pi/l_\mathrm{n}$ is element $\mathrm{n}$'s fundamental wave number. The first derivatives of the basis functions are
\begin{equation}
    \partial_qb^{(\text{FR})}_{\mathrm{n},\beta} = \beta k_\mathrm{n}
    \begin{cases}
    \cos\left(\beta k_\mathrm{n}q'\right)&; \beta > 0\\
    -\sin\left(\beta k_\mathrm{n}q'\right)&; \beta \leq 0
    \end{cases}.
\end{equation}
Then, we directly matching terms in the derivatives to the basis functions as,
\begin{align}
    \partial_qb^{(\text{FR})}_{\mathrm{n},\beta} &= \beta k_\mathrm{n}
    \begin{cases}
    \cos\left(\beta k_\mathrm{n}q'\right)&; -\beta < 0\\
    -\sin\left(\beta k_\mathrm{n}q'\right)&; -\beta \geq 0
    \end{cases}\nonumber\\
    &= \beta k_\mathrm{n}
    \begin{cases}
    \cos\left(-\beta k_\mathrm{n}q'\right)&; -\beta \leq 0\\
    \sin\left(-\beta k_\mathrm{n}q'\right)&; -\beta > 0
    \end{cases}\\
    &= \beta k_\mathrm{n}b^{(\text{FR})}_{\mathrm{n},-\beta}\,.
\end{align}
Then, we can construct differentiation matrix from this relation. On another hand, one can also directly apply \cref{eq:derivative_projection} to achieve the same result. 

Similarly, one can utilize trigonometric product properties as
\begin{align}
    b^{(\text{FR})}_{\mathrm{n},\beta_1}b^{(\text{FR})}_{\mathrm{n},\beta_2} &=
    \begin{cases}
    s_1s_2&; \beta_1 > 0,\beta_2 > 0\\
    s_1c_2&; \beta_1 > 0,\beta_2 \leq 0\\
    c_1s_2&; \beta_1 \leq 0,\beta_2 > 0\\
    c_1c_2&; \beta_1 \leq 0,\beta_2 \leq 0
    \end{cases}\nonumber\\
    &=\dfrac{1}{2}
    \begin{cases}
    c_--c_+&; \beta_1 > 0,\beta_2 > 0\\
    s_++s_-&; \beta_1 > 0,\beta_2 \leq 0\\
    s_+-s_-&; \beta_1 \leq 0,\beta_2 > 0\\
    c_-+c_+&; \beta_1 \leq 0,\beta_2 \leq 0
    \end{cases}
\end{align}
where $s_\square=\sin(\beta_\square k_\mathrm{n}q')$, $c_\square=\cos(\beta_\square k_\mathrm{n}q')$, and$\beta_\pm=\beta_1\pm\beta_2$. Then, matching each term according to
\begin{align}
    s_\square &= 
    \begin{cases}
        b^{(\text{FR})}_{\mathrm{n},\beta_\square}&; \beta_\square > 0\\
        b^{(\text{FR})}_{\mathrm{n},-\beta_\square}&; \beta_\square < 0
    \end{cases}\nonumber\\
    c_\square &= 
    \begin{cases}
        b^{(\text{FR})}_{\mathrm{n},\beta_\square}&; \beta_\square \leq 0\\
        b^{(\text{FR})}_{\mathrm{n},-\beta_\square}&; \beta_\square \geq 0
    \end{cases}\,.
\end{align}
Then, we can construct multiplication tensor from this relation. Like in the derivative case, one can also directly apply \cref{eq:multiplication_projection} to get the same result. However, in practice, the tensor would be truncated by maximum basis order. Since the element size and wavelength are independent, the basis functions are not orthogonal. This means that the term with higher $\abs{\beta}$ than the cut-off can have non-trivial contribution if one directly applying \cref{eq:multiplication_projection} causing instability. Hence, it is better to use the derived close form.

The boundary projection tensor is quite straight forward since the boundary for one-dimensional element is just a scalar, i.e., only one basis function such as delta function $\delta(q-q_{\mathrm{n}\circ})$, where $q_{\mathrm{n}\circ}$ is the coordinate of the boundary between $\Omega_\mathrm{n}$ and $\Omega_\circ$, is needed. This means that \cref{eq:boundary_projection} can be simplified to just an evaluation at the boundary
\begin{equation}
    B^{(\text{FR})}_{\mathrm{n}\circ,\beta} = \left.b^{(\text{FR})}_{\mathrm{n},\beta}\right|_{q=q_{\mathrm{n}\circ}}\,.
\end{equation}

One can directly apply a similar derivation to higher-dimensional cases, e.g., Fourier for each dimension. However, the calculation may be less straightforward, although a closed-form expression can still be obtained.

\subsubsection{Lagrange Polynomial Basis}
Lagrange polynomials are defined by their zeroes. Let the zeroes of a basis of one-dimensional element $\Omega_\mathrm{n}$ be at distinct points $q_{\mathrm{n}, 1}$, $\cdots$, $q_{\mathrm{n}, B}$. In this case, the basis functions with $\Omega_\mathrm{n}$ support for element $\mathrm{n}$ are defined as
\begin{equation}
    b^{(\text{LP})}_{\mathrm{n},\beta}(q) = \prod_{\beta'\neq\beta}\dfrac{q-q_{\mathrm{n}, \beta'}}{q_{\mathrm{n}, \beta}-q_{\mathrm{n}, \beta'}}.
\end{equation}
Similar to what we did for Fourier basis, one can try taking derivatives and rearranging terms such that the derivatives can be written as a linear combination of the basis functions. The derivatives of $b_{\mathrm{n},\beta}$ has the form of
\begin{equation}
    \partial_qb^{(\text{LP})}_{\mathrm{n},\beta} = b^{(\text{LP})}_{\mathrm{n},\beta}\sum_{\beta'\neq\beta}\dfrac{1}{q-q_{n,\beta'}}\,
\end{equation}
where we abuse the division by zero with $q-q_{n,\beta'}$ as long as it cancel the same term in $b_{\mathrm{n},\beta}^{(\text{LP})}$. Since the zeroes are distinct, each basis function is linearly independent. Together with the fact that derivative only reduce polynomial order, the term matching is guaranteed. 

In practice, we can choose Gauss-Lobatto-Legendre (GLL) points as the zeros to utilize Gaussian quadrature rule for faster integration with perfect accuracy up to $2B-3$ polynomial order. It turns out that, with this choice of zeros, Lagrange polynomial basis is also an interpolation basis which means that the projection of any function (with polynomial order at most $2B-3$) to each basis function is equal to the evaluation of that function at the non-zero point of that basis function~\cite{canuto2006spectral,deville2003high}, i.e.,
\begin{equation}
    \partial_qb^{(\text{LP})}_{\mathrm{n},\beta} = \sum_{\beta'}\left.\partial_qb^{(\text{LP})}_{\mathrm{n},\beta}\right|_{q=q_{\mathrm{n},\beta'}}b^{(\text{LP})}_{\mathrm{n},\beta'}\,.
\end{equation}

Similarly, the multiplication tensor can be achieve by evaluating the products between basis functions at different zeroes points
\begin{equation}
    b^{(\text{LP})}_{\mathrm{n},\beta_{1}}b^{(\text{LP})}_{\mathrm{n},\beta_{2}} \approx \sum_{\beta'}\left.b^{(\text{LP})}_{\mathrm{n},\beta_{1}}\right|_{q=q_{\mathrm{n},\beta'}}\left.b^{(\text{LP})}_{\mathrm{n},\beta_{2}}\right|_{q=q_{\mathrm{n},\beta'}}b^{(\text{LP})}_{\mathrm{n},\beta'}\,.
\end{equation}

As for the boundary projection, the result is exactly the same form as Fourier basis, i.e.,
\begin{equation}
    B^{(\text{LP})}_{\mathrm{n}\circ,\beta} = \left.b^{(\text{LP})}_{\mathrm{n},\beta}\right|_{q=q_{\mathrm{n}\circ}}\,.
\end{equation}

\subsubsection{Other Bases}
Fourier and Lagrange polynomial bases are special cases where derivative, multiplication and boundary projection tensors have their closed form. However, this doesn't mean that other bases will not work. They just need to be properly projected to the original basis, i.e., follow the exact calculation of \cref{eq:derivative_projection,eq:multiplication_projection,eq:boundary_projection}, and treated for any projection instability.

\section{Unified Residual Construction}\label{sec:unify_residual}
Since the solving scheme is a variational problem, any additional requirements for the solution can be integrated into the optimization target by framing the requirement as another variational problem to be added to the target function or loss function.

\subsection{Inter-element Continuity Residual}
With the boundary projection tensor, we can penalize the solution with how much discontinuity it has. This is the point when we can utilize the boundary projection tensor we developed. Basically, the discontinuity can be measure with how much the solution differs at the boundary between two adjacent elements. In other words, if element $\mathrm{n}_1$ and $\mathrm{n}_2$ are neighboring elements, the inter-element continuity residual is the difference between the boundary projection of both elements to their common boundary. 

We can also extend to the concept toward their other derivatives as well as shown in \cref{fig:bc}a where each element got differentiated for $\mathrm{p}$ times, i.e., the diagram is for the continuity of the $\mathrm{p}^\text{th}$ derivative.

\subsection{Boundary Condition Residuals}
The boundary condition can be viewed as the continuity of the solution to the environment element, i.e., complement of the domain. In most case, this is supplied by the problem. Hence, the boundary condition residual is exactly the same as the inter-element continuity case with one branch substituted by the environment boundary value as shown in \cref{fig:bc}b where the element $\mathrm{n}$'s $\mathrm{p}^\text{th}$ derivative is conditioned by the environment $\Omega_0$. 

\subsection{Combined Variational Objective}
\begin{figure}[htb!]
    \centering
    \includegraphics[width=\columnwidth]{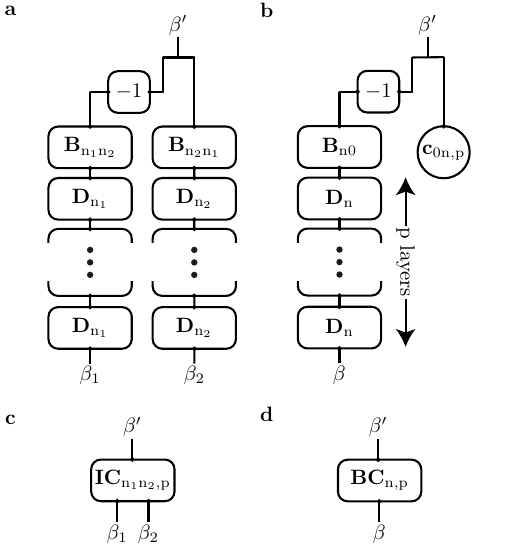}
    \caption{TNs representing element boundary constraint residual. \textbf{(a)} Inter-element continuity residual of the $\mathrm{p}^\text{th}$ derivative between elements $\mathrm{n}_1$ and $\mathrm{n}_2$. \textbf{(b)} Boundary condition of the $\mathrm{p}^\text{th}$ derivative of element $\mathrm{n}$ at the boundary with environment (element $0$) where $\mathbf{c}_{0\mathrm{n},\mathrm{p}}$ is the $\mathrm{p}^\text{th}$ derivative coefficient tensor of the boundary between $\Omega_\mathrm{n}$ and $\Omega_\mathrm{n}$ required by the problem. \textbf{(c)} and \textbf{(d)} are building block tensor nodes representing \textbf{a} and \textbf{b}, respectively.}
    \label{fig:bc}
\end{figure}

The complete variational objective is the sum of the weighted residual norm of the equation (PDE), inter-element continuities (IC), and boundary conditions (BC) where the norms need to be applied individually for all conditions. In other words, if we define the networks in \cref{fig:bc}a-b as in \cref{fig:bc}c-d, respectively, the final variational target is
\begin{equation}\label{eq:final_target}
    \mathcal{J}(\mathbf{c}) = \mathcal{J}^{(\mathrm{PDE})}(\mathbf{c}) + \mathcal{J}^{(\mathrm{IC})}(\mathbf{c}) + \mathcal{J}^{(\mathrm{BC})}(\mathbf{c})
\end{equation}
where
\begin{align}
    \mathcal{J}^{(\mathrm{PDE})} &= W^{(\mathrm{PDE})}\left(\mathbf{h}^{(\mathrm{PDE})\dagger}\odot \mathbf{h}^{(\mathrm{PDE})}\right)\nonumber\\
    \mathcal{J}^{(\mathrm{IC})} &= \sum_{\mathrm{n}_1,\mathrm{n}_2,\mathrm{p}} W^{(\mathrm{IC})}_{\mathrm{n}_1\mathrm{n}_2,\mathrm{p}}\left(\mathbf{h}_{\mathrm{n}_1\mathrm{n}_2,\mathrm{p}}^{(\mathrm{IC})\dagger}\odot \mathbf{h}^{(\mathrm{IC})}_{\mathrm{n}_1\mathrm{n}_2,\mathrm{p}}\right)\\
    \mathcal{J}^{(\mathrm{BC})} &= \sum_{\mathrm{n},\mathrm{p}} W^{(\mathrm{BC})}_{\mathrm{n},\mathrm{p}}\left(\mathbf{h}_{\mathrm{n},\mathrm{p}}^{(\mathrm{BC})\dagger}\odot \mathbf{h}^{(\mathrm{BC})}_{\mathrm{n},\mathrm{p}}\right)\,.\nonumber
\end{align}
\begin{figure*}[htb!]
    \centering
    \includegraphics[width=\textwidth]{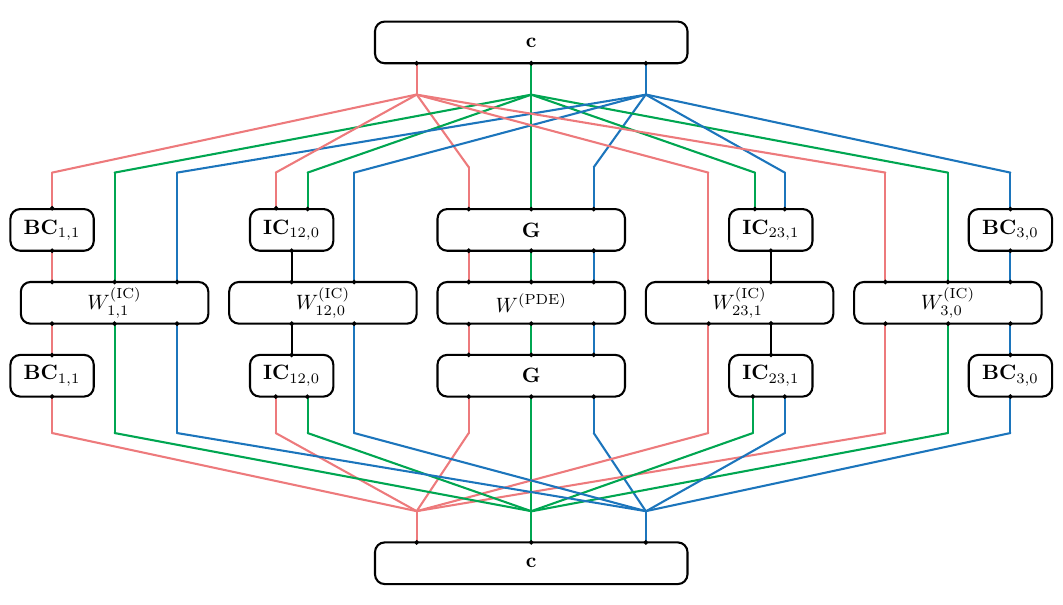}
    \caption{TN representing variational target of a system with 3 elements. Each element contraction is represented with different color bond: 1 (red), 2 (green), 3 (blue). The network composes of 5 branches each corresponding to different weighted residual norm. From left to right, the branches represent boundary condition on first derivative of element 1, continuity between element 1 and 2, main PDE, continuity of the first derivative between element 2 and 3, and boundary condition on element 3.}
    \label{fig:full_target}
\end{figure*}
Here, $\mathbf{h}$'s are residual terms 
\begin{align}
    \mathbf{h}^{(\mathrm{PDE})}(\mathbf{c}) &= \mathbf{G}\odot\mathbf{c}\nonumber\\
    \mathbf{h}^{(\mathrm{IC})}_{\mathrm{n}_1\mathrm{n}_2,\mathrm{p}}(\mathbf{c}) &= \mathbf{IC}_{\mathrm{n}_1\mathrm{n}_2,\mathrm{p}}\odot\mathbf{c}\\
    \mathbf{h}^{(\mathrm{BC})}_{\mathrm{n},\mathrm{p}}(\mathbf{c}) &= \mathbf{BC}_{\mathrm{n},\mathrm{p}}\odot\mathbf{c}\nonumber
\end{align}
and $W$'s are positive-definite metrics of  corresponding residuals. We also drop out some tensor indices for cleanliness as seen from the definition of $\mathbf{h}^{(\mathrm{PDE})}$ and $W^{(\mathrm{PDE})}$ compare to \cref{eq:aoe_res,eq:aoe_wrn}. \Cref{fig:full_target} illustrates an example TN that completely captures the variational target of a 3-element system with some continuities and boundary conditions.

\section{Target Optimization}\label{sec:target_optimization}
Given that the solution that satisfies all continuities and boundary conditions exist, minimizing \cref{eq:final_target} would give that solution. However, with the size complexity of the Fock space, the direct optimization is not viable. Hence, the application of TN ansatz is proposed to handle the data structure. This has a consequence that the optimizing space is no longer convex for weighted residual norm. Solving this problem is an active research in TN where we will directly use their results here. For the purpose of later demonstration, we will focus on one-dimensional problems where MPS ansatz is well studied and found success in solving quantum mechanical problem with DMRG algorithm.

\subsection{MPS Representation}
As shown in \cref{fig:tn_ansatz}b, MPS ansatz contains key features of one-dimensional problem with tensor nodes representing our finite elements align in a line with bonds connecting only adjacent nodes. With this ansatz, the Fock space coefficient can be approximated as
\begin{equation}
    \mathbf{c}^{\left[\beta_{\mathrm{n}\mathrm{j}}\right]}\approx\mathbf{c}_{1,\sigma_{12}}^{\left[\beta_{1\mathrm{j}}\right]}\mathbf{c}_{2,\sigma_{12}\sigma_{23}}^{\left[\beta_{2\mathrm{j}}\right]}\ldots\mathbf{c}_{N,\sigma_{N-1,N}}^{\left[\beta_{N\mathrm{j}}\right]}\,.
\end{equation}
The physical indices $\beta_{\mathrm{n}\mathrm{j}}$'s label the local coefficients associated with element $\mathrm{n}$, including any empty slots required by the selected maximum correlation order $J$. The bond dimension $S$ controls the number of correlations retained between the adjacent elements in the chain.

\subsection{Effective Single-site Problem}
The key concept of DMRG algorithm is to perform multiple single-site optimizations where all but one of the tensor node $\mathbf{c}_\mathrm{n}^{\left[\beta_{\mathrm{n}\mathrm{j}}\right]}$ are fixed. If we group up all $W$'s, $\mathbf{G}$, $\mathbf{IC}$'s, $\mathbf{BC}$'s, and $\mathbf{F}$'s in \cref{eq:final_target} as effective operator tensor $\mathbf{H}$, then the variational target is 
\begin{equation}\label{eq:sim_full_target}
    \mathcal{J}(\mathbf{c}) = \mathbf{c}^{\dagger\left[\beta'_{\mathrm{n}\mathrm{j}}\right]}\mathbf{H}_{\left[\beta'_{\mathrm{n}\mathrm{j}}\right],\left[\beta_{\mathrm{n}\mathrm{j}}\right]}\mathbf{c}^{\left[\beta_{\mathrm{n}\mathrm{j}}\right]}\,.
\end{equation}
During a single-site update, contracting the fixed part of the network produces an effective quadratic objective
\begin{equation}
    \mathcal{J}_{\text{eff},\mathrm{n}}\left(\mathbf{c}_\mathrm{n}\right) = \mathbf{c}_\mathrm{n}^\dagger\odot\mathbf{H}_{\text{eff},\mathrm{n}}\odot\mathbf{c}_\mathrm{n}\,,
\end{equation}
where $\mathbf{H}_{\text{eff},\mathrm{n}}$ is the contraction of $\mathbf{H}$ with all fixed nodes and the remaining contractions with $\mathbf{c}_\mathrm{n}$ are implicitly assumed in $\odot$ for cleanliness. With this local convexity, the update may be solved as either a local least-squares problem or explicit ``matrix'' equation
\begin{equation}
    \mathbf{H}_{\text{eff},\mathrm{n}}\odot\mathbf{c}_\mathrm{n} = \mathbf{0}\,.
\end{equation}

\subsection{Environment contractions}
\begin{figure}[htb!]
    \centering
    \includegraphics[width=\columnwidth]{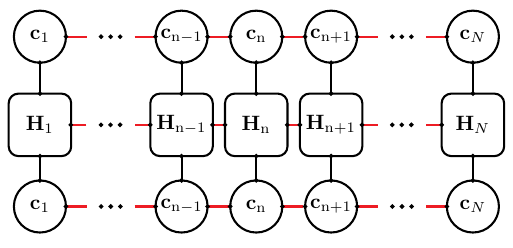}
    \caption{Tensor Network representation of full variational target (\cref{eq:sim_full_target}) under MPS ansatz for both Fock space coefficient tensor $\mathbf{c}$ and effective operator tensor $H$. }
    \label{fig:mps_target}
\end{figure}

Since the chosen TN ansatz stems from the nature of the expected solution, the PDE, continuities, and boundary conditions should support the same nature, i.e., $\mathbf{H}$ should have the same or similar structure as the chosen ansatz. This means that \cref{eq:sim_full_target} can be represented with TN similar to the one shown in \cref{fig:mps_target}.

We can obtain $\mathbf{H}_{\text{eff},\mathrm{n}}$ by contracting all fixed tensors on the left and right of the active site $\mathbf{c}_\mathrm{n}$ into $\mathbf{H}$. These contractions can be initialized recursively and updated dynamically through out the sequence of single-site updates to avoids redundancy of reconstructing the full contraction from scratch for every update.

\subsection{DMRG algorithm}
The Density Matrix Renormalization Group (DMRG) is a variational tensor-network algorithm initially developed for finding low-energy states, especially the ground state, of one-dimensional many-body quantum systems where the target wave function is represented as an MPS. The tensors are optimized sequentially back and forth along the chain. Since our target is also a variational optimization with MPS ansatz, the DMRG algorithm shown in \cref{alg:dmrg} should work for our case as well.
\begin{figure}[htb!]
\begin{algorithm}[H]
    \caption{DMRG}
    \label{alg:dmrg}
    \begin{algorithmic}
        \State Initialize $\mathbf{c}$
        \For {$\mathrm{n} = 1$ \textbf{to} $N$}
            \State $\mathbf{E}_{\text{L},\mathrm{n}} \gets$ Contract TN from element $\mathrm{n}-1$ to the left
            \State $\mathbf{E}_{\text{R},\mathrm{n}} \gets$ Contract TN from element $\mathrm{n}+1$ to the right
        \EndFor
        \For {$\mathrm{r} = 1$ \textbf{to} $R$}
            \For {$\mathrm{n} = 1$ \textbf{to} $N$}
                \State $\mathbf{H}_{\text{eff},\mathrm{n}} \gets$ Contract $\mathbf{H}_\mathrm{n}$, $\mathbf{E}_{\text{L},\mathrm{n}}$ and $\mathbf{E}_{\text{R},\mathrm{n}}$ together
                \State $\mathbf{c}_\mathrm{n} \gets \arg\min_{\mathbf{c}'_\mathrm{n}}\mathbf{c}_\mathrm{n}^{\prime\dagger}\odot\mathbf{H}_{\text{eff},\mathrm{n}}\odot\mathbf{c}'_\mathrm{n}$
                \State $\mathbf{E}_{\text{L},\mathrm{n}+1}\gets$ Contract $\mathbf{E}_{\text{L},\mathrm{n}}$ with $\mathbf{c}_\mathrm{n}^{\dagger}\odot\mathbf{H}_\mathrm{n}\odot\mathbf{c}_\mathrm{n}$
                
            \EndFor
            \For {$\mathrm{n} = N$ \textbf{to} $1$}
                \State $\mathbf{H}_{\text{eff},\mathrm{n}} \gets$ Contract $\mathbf{H}_\mathrm{n}$, $\mathbf{E}_{\text{L},\mathrm{n}}$ and $\mathbf{E}_{\text{R},\mathrm{n}}$ together
                \State $\mathbf{c}_\mathrm{n} \gets \arg\min_{\mathbf{c}'_\mathrm{n}}\mathbf{c}_\mathrm{n}^{\prime\dagger}\odot\mathbf{H}_{\text{eff},\mathrm{n}}\odot\mathbf{c}'_\mathrm{n}$
                \State $\mathbf{E}_{\text{R},\mathrm{n}-1}\gets$ Contract $\mathbf{E}_{\text{R},\mathrm{n}}$ with $\mathbf{c}_\mathrm{n}^{\dagger}\odot\mathbf{H}_\mathrm{n}\odot\mathbf{c}_\mathrm{n}$
            \EndFor
        \EndFor
        \State \Return $\mathbf{c}$
    \end{algorithmic}
\end{algorithm}
\end{figure}

Given effective operator tensor $\mathbf{H}$ and number of optimization round $R$, we start by initialize the Fock space coefficient tensor $\mathbf{c}$. Then, we prepare the left and right contractions of the TN for each element $\mathrm{n}$ denoted by $\mathbf{E}_{\text{L},\mathrm{n}}$ and $\mathbf{E}_{\text{R},\mathrm{n}}$, respectively. The preparation can be done recursively with dynamic memory. Then, we perform $R$ full loop sweep (back and forth) where at each step the single-site optimization is perform and a left or right contraction is updated for the next step. The $R$-loop sweeping can be replaced by stopping criteria.

\section{Complexity Analysis}\label{sec:complexity}
\subsection{Parameter Count}\label{subsec:space_complexity}
For our setup with $N$ elements, each with at most $2$ neighbors, $B$ elemental basis size, $S$ maximum bond dimension (for both $\mathbf{c}$ and $\mathbf{H}$) and $J$ maximum correlation order, the storage requirement are
\begin{enumerate}
    \item $\mathbf{c}$: $\mathcal{O}\left(NB^JS^2\right)$
    \item $\mathbf{H}$: $\mathcal{O}\left(NB^{2J}S^2\right)$
    \item $\mathbf{E}$'s: $\mathcal{O}\left(2NS^3\right)$
\end{enumerate}
for the DMRG process. For the preparation of $\mathbf{H}$, we can use the same storage as $\mathbf{H}$ by slowly contract each corresponding term of $\mathbf{L}$, continuities, and boundary conditions via summation filters. This gives the total of $\mathcal{O}\left(NS^2\left(B^{2J}+S\right)\right)$ storage space complexity.

\subsection{Runtime Complexity}
With the same parameters as in \cref{subsec:space_complexity}, $C$ conditions at each element edge, and all derivative order matrices are pre-calculated, there would be $\mathcal{O}\left(3^JNJ\right)$ terms in $\mathbf{L}$ due to the MPS ansatz. Each summation takes 2 matrices with size at most the same as $\mathbf{H}$ and add them through filter contraction to get an output of size at most the same as $\mathbf{H}$. The time complexity of $\mathbf{H}$ preparation is $\mathcal{O}\left(\left(C+3^JNJ\right)N^3B^{6J}S^6\right)$ where the first term is from continuities and boundary conditions while the second term if from PDE. This bound can be improve if the structure of $\mathbf{L}$ is known. For example, in the case of PDE problem, $\mathbf{L}$ only correlated with in the same elements, i.e., the number of possible terms reduces to $NJ$ and for each element the terms can be added in local space reducing the run time for the PDE part to $\mathcal{O}\left(NJB^{6J}S^6\right)$. The locality of continuities and boundary conditions would reduce the first term as well.

For the DMRG, the preparation of $\mathbf{E}$'s is iteratively contract a $\mathcal{O}(S^3)$-size, two $\mathcal{O}(B^JS^2)$-size, and a $\mathcal{O}(B^{2J}S^2)$-size tensors for at most $2N$ times with the optimal contraction time of $\mathcal{O}(NB^{2J}S^4)$. After that, there are at most $2NR$ single-site optimization each with the run time of $\mathcal{O}\left(B^{2J}S^6\right)$, $\mathcal{O}\left(B^{J\omega}S^{2\omega}\right)$, and $\mathcal{O}\left(B^{2J}S^4\right)$ for each sub-step in \cref{alg:dmrg} with $\omega\approx 2.37$ for the fastest matrix inversion algorithm~\cite{pan2017fast}. Therefore, the total run time of our method is $\mathcal{O}\left(\left(\left(C+3^JNJ\right)N^3B^{6J}+NRB^{\omega J}\right)S^6\right)$.

\section{Demonstration}\label{sec:demonstration}
\begin{figure}[!b]
    \centering
    \includegraphics[width=\columnwidth]{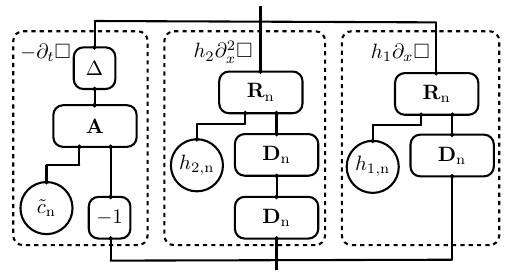}
    \caption{Tensor circuit for \cref{eq:space-dependent diffusion} where the circuit only show part of $\mathbf{G}$ from the residual (\cref{eq:aoe_res}) that act on element $\mathrm{n}$ tensor node of the Fock space coefficient tensor. $\mathbf{h}_1$, $\mathbf{h}_2$, $\widetilde{\mathbf{c}}$ are $-\sin(x)/2$, $(1+\cos(x)/2)$, and solution at the previous time step in Fock space representation. Note that we omit the network edges that connect these terms to the adjacent element for cleanliness.}
    \label{fig:diffuse_tc}
\end{figure}

We'll demonstrate the procedure of converting IVP to TN with the following diffusion equation with space-dependent diffusivity
\begin{equation}\label{eq:space-dependent diffusion}
    \partial_tu = \partial_x((1+\cos(x)/2)\partial_xu)
\end{equation}
with $u = u(t, x)$, $\Omega = [0, 2\pi]$, $u(0, x)=\cos(5x)/10$. We can rewrite the IVP in our format (\cref{eq:ivp_form}) as
\begin{equation}\label{eq:nlIVP}
    \dfrac{1}{\Delta t}(u^{(\mathrm{t})}-u^{(\mathrm{t}-1)}) = \left(1+\dfrac{\cos(x)}{2}\right)\partial^2_xu^{(\mathrm{t})} - \dfrac{\sin(x)}{2}\partial_xu^{(\mathrm{t})}\,.
\end{equation}
The corresponding tensor circuit of this IVP is as shown in \cref{fig:diffuse_tc}.

Then, we solve the IVP using $N$ uniform finite elements on $\Omega$, each equipped with $B$ local Lagrange polynomial basis functions, and an MPS ansatz with bond dimension $S$ for the coefficient tensor. The time evolution is advanced using the implicit Euler scheme described with time step $\Delta t$. At each time step, the DMRG sweep algorithm minimizes the weighted residual $\mathcal{J}(\mathbf{c})$ over the MPS variational manifold.

\begin{figure*}[!htb]
    \centering
    \includegraphics[width=0.9\textwidth]{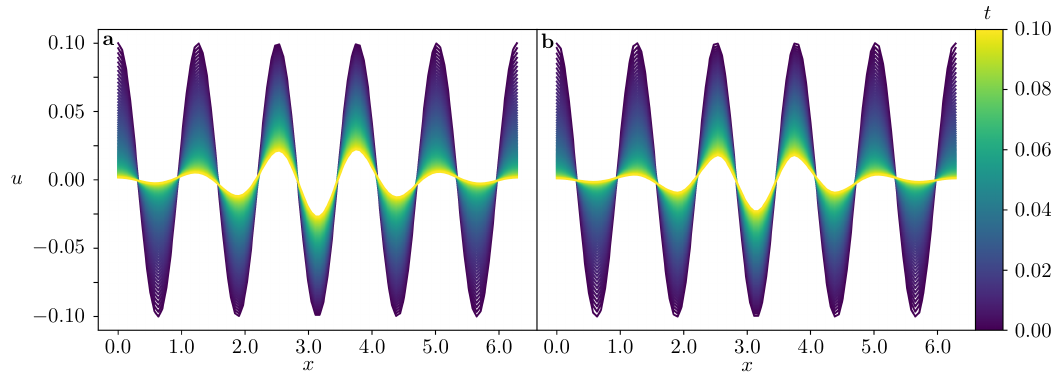}
    \caption{Numerical simulation results for the space-dependent diffusion equation $\partial_t u=\partial_x[\kappa(x)\partial_x u]$ with $\kappa(x)=1+\tfrac{1}{2}\sin x$, periodic domain $[0,2\pi]$, $\Delta t=0.001$ solving by \textbf{(a)} standard explicit Runge–Kutta method of order 5(4)~\cite{DORMAND198019} and \textbf{(b)} TN-FE variational method.}
    \label{fig:diffuse}
\end{figure*}

\begin{figure}[t!]
    \centering
    \includegraphics[width=\columnwidth]{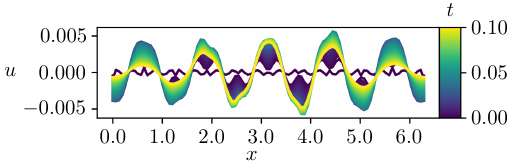}
    \caption{Discrepancy between the results from standard explicit Runge–Kutta method of order 5(4) and TN-FE variational method from \cref{fig:diffuse}.}
    \label{fig:be}
\end{figure}

\Cref{fig:diffuse} summarize the time-evolution results for explicit Runge-Kutta IVP simulation~\cite{DORMAND198019} (a) and TN variational method (b). The TN variational method with $N=10$, $B=4$, $S=1$ produce a reliable results with the difference of less than five percent of the original amplitude of the sinusoidal as shown in \cref{fig:be}. We also see that the continuity and smoothness of the results are also respected to within numerical precision confirming that the variational target encodes the boundary constraints.

\subsection{Convergence and Accuracy}
\begin{figure}[!htb]
    \centering
    \includegraphics[width=\columnwidth]{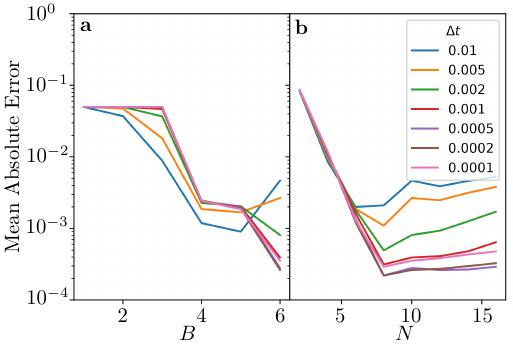}
    \caption{Mean absolute error trend of TN variational method for solving space-dependent diffusion IVP with different \textbf{(a)} basis order $B$ and \textbf{(b)} number of elements $N$ at different time step size $\Delta t$. For each fixed $\Delta t$, the trends exhibit convergence up to a level after which the better temporal resolution is required for a more accurate result.}
    \label{fig:convergence}
\end{figure}

To assess the convergence properties of the TN solver, we study the behavior of the solution error as functions of Lagrange polynomial basis order $B$, and the number of finite elements $N$ at different time step size $\Delta t$. Since both $B$ and $N$ are, in some way, representation of spatial resolution, we should expect the behavior of the trend to be similar to standard finite difference method which exhibit the stability condition where the convergence cannot be improved purely by increasing spatial resolution while the temporal resolution is fixed. This is clearly shown in \cref{fig:convergence}.

Together, these results demonstrate that the TN-FE framework introduced in the companion paper can produce accurate IVP solution with the computational scaling controlled by the bond dimension rather than the full exponential size of the coefficient space. While the examples here are restricted to 1D MPS/DMRG for clarity, the formalism is directly extensible to higher-dimensional problems via PEPS or tree tensor network ansatzes.

\section{Conclusion}\label{sec:conclusion}

We have developed a practical tensor-network finite-element framework for solving PDEs with polynomial computational scaling under bounded local correlation order and tensor-network bond dimension. Building on the general TN-FE formulation introduced in our companion work, we reorganize the ingredients of conventional finite-element discretization into a common variational construction. Rather than requiring the elemental basis functions to assemble into a globally conforming approximation space, we treat the coefficients of different elements as independent variational degrees of freedom and impose the governing PDE, inter-element continuity, and boundary conditions through separate residual operators within a single weighted objective. Spatial derivatives, multiplication, constant factors, and boundary projections are represented as reusable local tensor operators, from which more complicated differential expressions can be constructed. We explicitly derive these operators for Fourier and Lagrange polynomial bases, while the same projection procedure allows other local approximation spaces to be incorporated without changing the overall optimization framework.

The exponentially large augmented coefficient space generated by nonlinear PDEs is made computationally accessible by restricting the coefficient tensor to a tensor-network manifold. For one-dimensional systems, an MPS representation reduces the storage requirement from exponential in the number of elements to polynomial scaling, while DMRG converts the resulting global non-convex optimization into a sequence of convex single-site subproblems. By exploiting locality and recursively updated environment contractions, the PDE, continuity, and boundary contributions can be optimized without explicitly constructing the full augmented coefficient tensor. Our complexity analysis shows that, for fixed correlation order, bond dimension, and local basis size, both storage and computational cost scale polynomially with the number of finite elements. The variable-coefficient diffusion example demonstrates that the resulting solver reproduces the expected time evolution with controlled error, preserves the imposed inter-element and boundary constraints, and exhibits systematic convergence under both element and basis refinement.

The polynomial scaling established here should not be interpreted as a universal computational advantage over conventional PDE solvers. Its practical efficiency depends on whether the solution and the corresponding variational operators admit compact tensor-network representations with sufficiently small bond dimensions and correlation order. The present implementation is also restricted to one-dimensional MPS/DMRG calculations, and the polynomial prefactors can become substantial as the local basis size, nonlinear correlation order, or required bond dimension increases. Future work should therefore investigate adaptive control of bond dimension and local approximation order, more efficient contraction and operator-assembly strategies, and extensions to higher-dimensional tensor-network geometries such as PEPS and tree tensor networks. Systematic benchmarks against established finite-element, spectral-element, and discontinuous-Galerkin solvers will also be important for identifying the classes of PDEs for which the correlated representation provides a practical computational advantage. Overall, the framework separates local approximation, physical constraints, and global correlation into distinct but compatible components of a single variational problem. This provides a flexible route from conventional finite-element discretization to tensor-network optimization and establishes an algorithmic foundation for polynomial-scaling, quantum-inspired, and potentially quantum-assisted PDE solvers.

\section*{Acknowledgements}
The authors thank Drs. Alan Edelman and Nuno Loureiro for insightful discussions. We dedicate this work to Dr. Loureiro's memory and note with sadness that he passed away during the course of this work. This work was supported by the National Science Foundation (NSF) Convergence Accelerator Award No. 2345084, the U.S. Department of Energy (DOE), Basic Energy Sciences (BES), Award No. DE-SC0020148, MIT Energy Initiative, and support from R. Wachnik.

\bibliographystyle{unsrt}
\bibliography{references}
\end{document}